\RequirePackage{ifpdf}
\ifpdf 
\documentclass[pdftex]{sigma}
\else
\documentclass{sigma}
\fi

\usepackage{euscript}
\usepackage{mathrsfs}
\usepackage{dsfont}
\usepackage{oldgerm}
\newfont{\Sans}{cmss10 scaled\magstep1}         \newcommand{\sans}[1]{\hbox{\Sans {#1}}}
\newcommand{\set}[1]{\left\{#1\right\}}

\def\dsr{{\mathds R}}

\def\dsz{{\mathds Z}}

\def\exp#1{{\rm exp}(#1)}

\font\tentt=cmbx12 scaled \magstephalf

\DeclareMathOperator{\spann}{\mathds{R}-span}
\DeclareMathOperator{\spannq}{\mathds{Q}-span}
\DeclareMathOperator{\zspann}{\mathds{Z}-span}
 \DeclareMathOperator{\ind}{Ind}
 \DeclareMathOperator{\Ad}{Ad}

\DeclareMathOperator{\bbase}{\mathscr{B}}
\DeclareMathOperator{\mult}{\mathbf{\tentt{m}}}
\DeclareMathOperator{\orbite}{\EuScript{O}}

\begin{document}


\renewcommand{\PaperNumber}{062}

\FirstPageHeading

\ShortArticleName{On the Moore Formula of Compact Nilmanifolds}

\ArticleName{On the Moore Formula of Compact Nilmanifolds}

\Author{Hatem HAMROUNI}

\AuthorNameForHeading{H. Hamrouni}

\Address{Department of Mathematics, Faculty of Sciences at Sfax,\\
Route  Soukra, B.P.~1171, 3000 Sfax, Tunisia}
\Email{\href{mailto:hatemhhamrouni@voila.fr}{hatemhhamrouni@voila.fr}}

\ArticleDates{Received December 17, 2008, in f\/inal form June 04, 2009;  Published online June 15, 2009}

\Abstract{Let $G$ be a connected and simply connected  two-step
nilpotent Lie group and
 $\Gamma$ a lattice  subgroup of $G$. In this note, we give a new
 multiplicity formula, according to the sense of Moore,
  of  irreducible unitary representations
 involved in the decomposition of the quasi-regular representation $\ind_\Gamma^G(1)$.
 Extending then the Abelian case.}


\Classification{22E27}

\section{Introduction}

 Let $G$ be a connected simply connected nilpotent Lie group
with Lie algebra $\frak g$ and suppose $G$ contains a discrete
cocompact subgroup $\Gamma$. Let $\sans{R}_\Gamma =
\ind_\Gamma^G(1)$ be the quasi-regular representation of $G$ induced
from $\Gamma$. Then $\sans{R}_\Gamma$ is direct sum of irreducible
unitary representations each occurring with f\/inite multiplicity~\cite{Gelfand1}; we will write
\begin{gather*}
\sans{R}_\Gamma =  \sum _{\pi \in (G:\Gamma)} \mult(\pi, G,
\Gamma, 1)
 \pi.
 \end{gather*}
A basic problem in representation theory is to determine the
spectrum  $(G: \Gamma)$ and the multiplicity function $\mult(\pi, G,
\Gamma, 1)$. C.C.~Moore f\/irst studied this problem in
\cite{Moore1}. More precisely, we have the following theorem.

\begin{theorem}\label{cond-Moore}
Let $G$ be a simply connected nilpotent Lie group with Lie algebra
$\frak g$ and  $\Gamma$  a lattice subgroup  of  $G$ $($i.e., $\Gamma$
is a discrete cocompact subgroup of $G$ and $\log(\Gamma)$ is an
additive subgroup of ${\frak g})$. Let $\pi$ be an irreducible unitary
representation with coadjoint orbit $\orbite_\pi^G$. Then $\pi$
belongs to $(G: \Gamma)$ if and only if $\orbite_\pi^G$ meets $
\frak g_\Gamma^* = \{l\in \frak g^*,\; \langle l, \log(\Gamma)\rangle
\subset \mathds{Z}\}$ where $\frak g^*$ denotes the dual space of
$\frak g$.
\end{theorem}

Later R.~Howe \cite{Howe2} and L.~Richardson \cite{Richardson1} gave
independently  the decomposition of $\sans{R}_\Gamma $ for an
arbitrary compact nilmanifold. In this paper, we pay attention to
the question wether the multiplicity formula
\begin{gather*}
\mult(\pi, G, \Gamma, 1)= \#[\orbite_\pi^G\cap \frak
g_\Gamma^*/\Gamma]\qquad \forall \, \pi \in (G:\Gamma)
\end{gather*}
required in the Abelian context, still holds for non commutative
nilpotent Lie groups (we write $\#A$ to denote the cardinal number
of a set $A$). In \cite{Moore1}, Moore showed the following
inequality
\begin{gather}\label{Moore-inequality}
\mult(\pi, G, \Gamma, 1)\leq \#[\orbite_\pi^G\cap \frak
g_\Gamma^*/\Gamma]\qquad \forall \, \pi \in (G:\Gamma),
\end{gather}
where $\Gamma$ is a lattice subgroup of $G$, and produced an example
for which the inequality (\ref{Moore-inequality}) is strict. More
precisely, he showed that
\begin{gather}\label{Moore-equality}
\mult(\pi, G, \Gamma, 1)^2= \#[\orbite_\pi^G\cap \frak
g_\Gamma^*/\Gamma]\qquad \forall \, \pi \in (G:\Gamma)
\end{gather}
in the case of the $3$-dimensional Heisenberg group and $\Gamma$ a
lattice subgroup.  The present paper aims to show that every
connected, simply connected two-step nilpotent Lie group satisf\/ies
equation (\ref{Moore-equality}). We present therefore a counter
example for $3$-step nilpotent Lie groups.

\section{Rational structures and uniform subgroups}

In this section, we summarize facts concerning rational structures
and uniform subgroups in a~connected, simply connected nilpotent Lie
groups. We recommend~\cite{Cor1} and~\cite{Onishchik2} as a~references.

\subsection{Rational structures}

 Let $G$ be a nilpotent, connected
and simply connected real Lie group and let $\frak g$ be its Lie
algebra. We say that $\frak g$ (or $G$) has a \textit{rational
structure} if there is a Lie algebra $\frak g_\mathds{Q}$ over
$\mathds{Q}$ such that $\frak g \cong \frak g_\mathds{Q} \otimes
\mathds{R}$. It is clear that $\frak g$ has a rational structure if
and only if $\frak g$ has an $ \mathds{R}$-basis
$\{X_1,\dots,X_n\}$ with rational structure constants.

 Let $\frak g$ have a f\/ixed
rational structure given by $\frak g_\mathds{Q}$ and let $\frak h$
be an $\mathds{R}$-subspace of~$\frak g$. Def\/ine $\frak
h_\mathds{Q}= \frak h \cap \frak g_\mathds{Q}$. We say that $\frak
h$ is \textit{rational} if $\frak h = \spann\set{\frak
h_\mathds{Q}}$, and that a connected, closed subgroup~$H$ of $G$ is
\textit{rational} if its Lie algebra $\frak h$ is rational. The
elements of $\frak g_\mathds{Q}$ (or $G_\mathds{Q}= \exp{\frak
g_\mathds{Q}}$) are called \textit{rational elements} (or
\textit{rational points}) of $\frak g$ (or $G$).\vskip0.3cm

\subsection{Uniform subgroups} A discrete subgroup $\Gamma$
is called \textit{uniform} in $G$ if the quotient space $G/\Gamma$
is compact. The homogeneous space $G/\Gamma$ is called a
\textit{compact nilmanifold}. A proof of the next result can be
found in Theorem 7 of \cite{Malcev1} or in Theorem 2.12 of
\cite{Raghunathan}.
\begin{theorem}[the Malcev rationality criterion]
Let $G$ be a simply connected nilpotent Lie group, and let $\frak g$
be its Lie algebra. Then $G$ admits a uniform subgroup $\Gamma$ if
and only if $\frak g$ admits a basis $\{X_1,\dotsc,X_n\}$ such that
\begin{gather*}
[X_i, X_j] =  \sum_{k=1}^n c_{ijk} X_k,\qquad \forall\,  1\leq i, j\leq n,
\end{gather*}
 where the constants $c_{ijk}$ are all
rational. $($The $c_{ijk}$ are called the structure constants of
$\frak g$ relative to the basis $\set{X_1, \ldots,
X_n}.)$
\end{theorem}

More precisely, we have, if $G$ has a uniform subgroup $\Gamma$,
then $\frak g$ (hence $G$) has a rational structure such that $\frak
g_\mathds{Q}= \spannq\set{\log(\Gamma)}$. Conversely, if $\frak g$
has a rational structure given by some $\mathds{Q}$-algebra $\frak
g_\mathds{Q}\subset \frak g$, then $G$ has a uniform subgroup
$\Gamma$ such that $\log(\Gamma) \subset \frak g_\mathds{Q}$ (see~\cite{Cor1} and~\cite{Malcev1}).  If we endow $G$ with the rational
structure induced by a~uniform subgroup $\Gamma$ and if $H$ is a Lie
subgroup of $G$, then $H$ is rational if and only if $H\cap \Gamma$
is a uniform subgroup of $H$. Note that the notion of rational
depends on $\Gamma$.

\subsection{Weak and strong Malcev basis}

 Let $\frak g$ be a
nilpotent Lie algebra and let $\bbase=\set{X_1, \ldots, X_n}$ be a
basis of $\frak g$.
 We say that $\bbase$ is a weak (resp. strong) Malcev
basis for $\frak g$ if  $\frak g_i=\spann\set{X_1, \ldots, X_i}$ is
a subalgebras (resp. an ideal) of $\frak g$ for each $1\leq i\leq n$
(see~\cite{Cor1}).

 Let $ \Gamma$ be a uniform subgroup of $G$. A
strong or weak Malcev  basis $\{X_1,\dotsc,X_n\}$ for $\frak g$ is
said to be \textit{strongly based on} $\Gamma$ if{\samepage
\begin{gather*}
\Gamma= \exp {\mathds{Z}X_1}\cdots \exp {\mathds{Z}X_n}.
\end{gather*}
 Such a basis always exists (see \cite{Malcev1,Cor1,Matsushima1}).}

A proof of the next result can be found in  Proposition~5.3.2 of
\cite{Cor1}.

\begin{proposition}
Let $\Gamma$ be uniform subgroup in a nilpotent Lie group $G$, and
let $H_1\subsetneqq H_2\subsetneqq\cdots\subsetneqq H_k =G$ be
rational Lie  subgroups of $G$. Let $\frak h_1, \ldots, \frak
h_{k-1}, \frak h_k =\frak g$ be the corresponding Lie algebras. Then
there exists a weak Malcev basis $\set{X_1,\ldots, X_n}$ for $\frak
g$ strongly based on $\Gamma$ and passing through $\frak h_1,
\ldots, \frak h_{k-1}$. If the $H_j$ are all normal, the basis can
be chosen to be a strong Malcev basis.
\end{proposition}

\subsection{Lattice subgroups}
\begin{definition}[\cite{Moore1}]\label{definition_lattice} Let $\Gamma$ be a uniform subgroup of a simply
connected nilpotent Lie group $G$, we say that $\Gamma$ is a lattice
subgroup of $G$ if  $ \log(\Gamma)$ is an Abelian subgroup of $\frak
g$.
\end{definition}

In \cite{Moore1}, Moore shows that if a simply connected nilpotent
Lie group $G$ satisf\/ies the Malcev rationality criterion, then $G$
admits a lattice subgroup.

 We close this section with the
following proposition~\cite[Lemma 3.9]{Cor9}.

\begin{proposition}\label{Imran} If $\Gamma$ is a lattice subgroup of a simply connected
nilpotent Lie group $G=\exp {\frak g}$ and $\set{X_1, \ldots, X_n}$
is a weak Malcev basis of $\frak g$ strongly based on $\Gamma$, then
$\set{X_1, \ldots, X_n}$ is a $\dsz$-basis for the additive lattice
$\log(\Gamma)$ in $\frak g$.
\end{proposition}

\section{Main result}

We begin with the following def\/inition.
\begin{definition} Let $G$ be a connected, simply connected nilpotent Lie
group which satisf\/ies the Malcev rationality criterion, and let
$\frak g$ be its Lie algebra.
\begin{itemize}\itemsep=0pt
\item[$(1)$]  We say that $G$ satisf\/ies the Moore formula at a lattice
subgroup $\Gamma$ if we have
\[
\mult(\pi, G, \Gamma, 1)^2 = \#[\orbite_\pi^G\cap \frak
g_\Gamma^*/\Gamma], \qquad \forall \, \pi \in (G:\Gamma)).
\]

\item[$(2)$]
 We say that $G$ satisf\/ies the Moore formula if $G$ satisf\/ies
  the Moore formula at  every lattice
subgroup $\Gamma$ of $G$.
\end{itemize}
\end{definition}

\noindent
\textbf{Examples.}
\begin{itemize}\itemsep=0pt
  \item[$(1)$] Every Abelian Lie group satisf\/ies the Moore formula.
  \item[$(2)$] The $3$-dimensional Heisenberg group satisf\/ies the
  Moore formula (see \cite[p.~155]{Moore1}).
\end{itemize}

The main result of this paper is the following theorem.

\begin{theorem}\label{Thm-Moore-Formula}
Every connected, simply connected two-step nilpotent Lie group
satisfies the Moore formula.
\end{theorem}

 Before proving Theorem~\ref{Thm-Moore-Formula}, we must review
more of the   Corwin--Greenleaf multiplicity formula.

\subsection[The Corwin-Greenleaf multiplicity formula]{The Corwin--Greenleaf multiplicity formula}

Using the Poisson summation and Selberg trace formulas, L.~Corwin
and F.P.~Greenleaf~\cite{Cor9} gave a formula for $\mult(\pi, G,
\Gamma, 1)$ that depended only on the coadjoint orbit in $\frak g^*$
corresponding to $\pi$ via Kirillov theory. We state their formula
for lattice subgroups. Let $\Gamma$ be a lattice subgroup of a
connected, simply connected nilpotent Lie group $G=\exp{\frak g}$.
Let
\[
\frak g_\Gamma^*=\set{l\in \frak g^*:\; \langle l,
\log(\Gamma)\rangle \subset \dsz}.
\]
Let $\pi_l$ be an irreducible
unitary representation of $G$ with coadjoint orbit
$\orbite_{\pi_l}^G\subset \frak g^*$ such that $\orbite_{\pi_l}^G\ne
\set{l}$. According to Theorem~\ref{cond-Moore}, we have
$\mult(\pi_l, G, \Gamma, 1)>0$ if and only if
$\orbite_{\pi_l}^G\cap\frak g_\Gamma^*\ne \emptyset$, so we will
suppose this intersection is nonempty. The set
$\orbite_{\pi_l}^G\cap\frak g_\Gamma^*$ is $\Gamma$-invariant. For
such $\Gamma$-orbit $\Omega\subset \orbite_{\pi_l}^G\cap\frak
g_\Gamma^*$ one can associate a number $c(\Omega)$ as follows:
 let $f\in \Omega$ and
$\frak g(f) = \ker(B_f)$, where $B_f$ is the skew-symmetric bilinear
form on $\frak g$ given by
\[
B_f(X, Y)=
\langle f, [X, Y]\rangle, \qquad X, Y \in \frak g.
\]
 Since $\langle f,
\log(\Gamma)\rangle\subset \dsz$ then $\frak g(f)$ is a rational
subalgebra. There exists a weak Malcev basis $\set{X_1, \ldots,
X_n}$ of $\frak g$ strongly based on $\Gamma$ and passing through
$\frak g(f)$ (see \cite[Proposition 5.3.2]{Cor1}). We write $\frak
g(f)= \spann\set{X_1, \ldots, X_s}$.
 Let
\begin{gather}\label{omi}
A_f = \mbox{Mat}\big(\langle f, [X_i, X_j]\rangle:\; s< i, j \leq
n\big). \end{gather}
Then $\det(A_f)$ is independent  of the basis
satisfying the above conditions and depends only on the
$\Gamma$-orbit $\Omega$. Set
\[
c(\Omega) = \big(\det(A_f)\big)^{-\frac{1}{2}}.
\]
Then $c(\Omega)$ is a positive rational number and the multiplicity
formula of  Corwin--Greenleaf is
\begin{gather}\label{mult-Corwin-Greenleaf}
\mult(\pi_l, G, \Gamma, 1) =\left\{
                             \begin{array}{ll}
                               1, & \hbox{if} \ \ \frak g(l) = \frak g, \\
                               \displaystyle  \sum\limits_{\Omega\in [\orbite_{\pi_l}^G\cap \frak g_\Gamma^*/\Gamma]}
c(\Omega), & \hbox{otherwise}.
                             \end{array}
                           \right. \end{gather}
 For details see~\cite{Cor9}.

\begin{proof}[Proof of Theorem~\ref{Thm-Moore-Formula}]
Let $l\in \orbite_\pi^G\cap \frak g_\Gamma^*$. The result is obvious
if $\frak g(l) = \frak g$. Next, we suppose that $\frak g(l) \ne
\frak g$. Since $G$ is two-step nilpotent Lie group then $\frak
g(l)$ is an ideal of $\frak g$, and hence we have $\frak g(l) =
\frak g(f)$ for every $f\in \orbite_\pi^G$ and $\orbite_\pi^G = l+
\frak g(l)^\bot$ (see \cite[Theorem 3.2.3]{Cor1}). On the other
hand, as $l$ belongs to $\frak g_\Gamma^*$ then $\frak g(l)$ is
rational. By Proposition~5.3.2 of \cite{Cor1} there exists a~Jordan--H\"{o}lder basis $\bbase=\set{X_1, \ldots, X_n}$ of $\frak g$
strongly based on $\Gamma$ and passing through $\frak g(l)$. Set
$\frak g(l)=\spann\set{X_1, \ldots, X_s}$.

Then, for every $\Omega \in [\orbite_\pi^G\cap \frak g_\Gamma^*/\Gamma]$ and for every $f\in \Omega$, we have
\begin{gather*}
c(\Omega) = \det(A_f)^{- \frac{1}{2}}=\det(A_l)^{- \frac{1}{2}}= c(\Gamma\cdot l),
\end{gather*}
since $f\vert_{[\frak g, \frak g]} = l\vert_{[\frak g, \frak g]}$.
It follows from (\ref{mult-Corwin-Greenleaf}) that
\begin{gather}\label{etiquette1}
\mult(\pi, G, \Gamma, 1) = \#[\orbite_\pi^G\cap \frak
g_\Gamma^*/\Gamma]\  c(\Gamma\cdot l).
\end{gather}

Next, we calculate $\#[\orbite_\pi^G\cap \frak g_\Gamma^*/\Gamma]$.
Let $(t_1, \ldots, t_n)\in \dsz^n$ and $f\in \orbite_\pi^G\cap \frak
g_\Gamma^*$. We have
\begin{gather*}
  \big(\exp{-t_1X_1}\cdots \exp{-t_nX_n}\big)\cdot f  =  f+ \sum_{i=s+1}^n
  \left(\sum_{j=s+1}^n t_j \langle f, [X_j, X_i]\rangle \right)X_i^*\\
\phantom{\big(\exp{-t_1X_1}\cdots \exp{-t_nX_n}\big)\cdot f }{}= f+ \sum_{i=s+1}^n
  \left(\sum_{j=s+1}^n t_j \langle l, [X_j, X_i]\rangle \right)
  X_i^*,
\end{gather*}
since $f\vert_{[\frak g, \frak g]} = l\vert_{[\frak g, \frak g]}$.
It follows that
\[
\Gamma\cdot f = f+ \sum_{j=s+1}^n \dsz e_j,
\]
where
\[
e_j=
  \sum_{i=s+1}^n  \langle l, [X_j, X_i]\rangle X_i^*, \qquad \forall\, s< j\leq n.
\]
Let
\[
\textswab{L} = \orbite_\pi^G\cap \frak
  g_\Gamma^*-f= \bigoplus_{s<i\leq n } \dsz X_i^*\qquad \mbox{and}\qquad
  \textswab{L}_0=  \sum_{j=s+1}^n \dsz e_j.
\]
  Since $\frak g(l) \cap \spann\set{X_{s+1}, \ldots, X_n}=\set{0}$, then
  the vectors  $e_{s+1}, \ldots, e_n$ are
linearly independent.
  Therefore,  $\textswab{L}_0$ is a sublattice of $\textswab{L}$. It is well known
  that there exist $ \varepsilon_{s+1}, \ldots,
\varepsilon_n$ a linearly independent vectors of $\frak g^*$ and
$d_{s+1}, \ldots, d_n \in \mathds{N}^*$ such that
\[
\textswab{L}= \bigoplus_{s<i\leq n } \dsz \varepsilon_i\qquad \mbox{and}\qquad \textswab{L}_0 =  \bigoplus_{s<i\leq n } d_i\dsz \varepsilon_i.
\]
Consequently, we have
\[
\#[\orbite_\pi^G\cap \frak g_\Gamma^*/\Gamma]= d_{s+1}\cdots d_n.
\]
Let  $[ \varepsilon_{s+1}, \ldots, \varepsilon_n]$ be the matrix
with column vectors $ \varepsilon_{s+1}, \ldots, \varepsilon_n$
expressed in the basis $(X_{s+1}^*$, $\ldots, X_n^*)$. From
\[
\textswab{L} = \bigoplus_{s<i\leq n } \dsz X_i^*=
 \bigoplus_{s<i\leq n } \dsz \varepsilon_i,
 \]
we deduce that
\[
[ \varepsilon_{s+1}, \ldots, \varepsilon_n]\in
\mathrm{GL}(n-s, \dsz).
\]
On the other hand, let $[e_{s+1}, \ldots,
e_n]$ (resp. $[ d_{s+1}\varepsilon_{s+1}, \ldots,
d_{n}\varepsilon_n]$) be the matrix with column vectors $e_{s+1},
\ldots, e_n$ (resp.\ $d_{s+1}\varepsilon_{s+1}, \ldots,
d_{n}\varepsilon_n$) expressed in the basis $(X_{s+1}^*, \ldots,
X_n^*)$. Since
\[
\textswab{L}_0=  \sum_{j=s+1}^n \dsz e_j =
 \bigoplus_{s<i\leq n } d_i\dsz \varepsilon_i,
 \]
  then there
exists $T\in \mathrm{GL}(n-s, \dsz)$ such that
\[
[e_{s+1}, \ldots,
e_n] = [ d_{s+1}\varepsilon_{s+1}, \ldots, d_{n}\varepsilon_n] T.
\]
 The latter condition can be written
\[
{}^tA_l = [\varepsilon_{s+1}, \ldots, \varepsilon_n]
\mathrm{diag}[d_{s+1}, \ldots, d_n] T .
\]
 Form this it follows that
\[
\det(A_l)=d_{s+1}\cdots d_n.
\]
Consequently
\begin{gather}\label{new}
\#[\orbite_\pi^G\cap \frak g_\Gamma^*/\Gamma]= \det(A_l).
\end{gather}
Substituting the last expression (\ref{new})  into
(\ref{etiquette1}), we obtain
\begin{gather*}
 \mult(\pi, G, \Gamma, 1)^2 =
\#[\orbite_\pi^G\cap \frak g_\Gamma^*/\Gamma].
\end{gather*}
This completes the proof.
\end{proof}

As a consequence of the above result, we obtain the
following result.
\begin{corollary}\label{description_Pesce} Let $G$ be a connected,
simply connected two-step nilpotent Lie group, let $\frak g$ be the
Lie algebra of $G$, and let $\Gamma$ be a lattice subgroup of $G$.
Let $l\in \frak g^*$ such that  the representation $\pi_l$ appears
in the decomposition of  $\sans{R}_\Gamma$. Let $A_l$ as in~\eqref{omi}. The multiplicity of $\pi_l$ is
\begin{gather*}\label{Pesce-formula}
\mult(\pi_l, G, \Gamma, 1)=\left\{
                             \begin{array}{ll}
                               1, & \hbox{if} \ \  \frak g(l) = \frak g, \\
                               (\det(A_l))^{\frac{1}{2}}, &
                               \hbox{otherwise}.
                             \end{array}
                           \right.
\end{gather*}
\end{corollary}

\begin{remark} Note that in  \cite{Pesce1}, H.~Pesce obtained the above
result more generally when $\Gamma$ is a uniform subgroup of $G$.
\end{remark}

\section{Three-step example}

In this section, we give an example of three-step nilpotent Lie
group that does not satisfy the Moore formula. Consider the
$4$-dimensional three-step nilpotent Lie algebra
\begin{gather*}
\frak g = \spann\set{X_1, \ldots, X_4}
\end{gather*}
with Lie brackets  given by
\begin{gather*}
[X_4, X_i] = X_{i-1}, \qquad i=2, 3,
\end{gather*}
and the non-def\/ined brackets being equal to zero or obtained by
antisymmetry. Let $G$ be the simply connected Lie group with Lie
algebra $\frak g$. The group $G$ is called the generic f\/iliform
nilpotent Lie group of dimension four. Let $\Gamma$ be the lattice
subgroup of $G$ def\/ined by
\begin{gather*}
\Gamma  =  \exp{\dsz X_1}\exp{\dsz X_2}\exp{\dsz X_3}\exp{6\dsz
X_4} =  \exp{\dsz X_1\oplus\dsz X_2\oplus\dsz X_3\oplus 6\dsz X_4}.
\end{gather*}
Let $l= X_1^*$. It is clear that the ideal $\frak m =
\spann\set{X_1, \ldots, X_3}$ is a rational polarization at $l$. On
the other hand, we have $\langle l, \frak m\cap \log(\Gamma)\rangle
\subset \dsz$. Consequently, the representation $\pi_l$ occurs in~$\sans{R}_\Gamma$ (see~\cite{Richardson1,Howe2}). Now, we
have to calculate $\#[\orbite_{\pi_l}^G\cap \frak
g_\Gamma^*/\Gamma]$.

Following \cite{Cor1} or \cite{Nielsen5}, the  coadjoint orbit of
$l$  has  the form
\[
\orbite_{\pi_l}^G = \set{X_1^*+  t X_2^* +
 \frac{t^2}{2}  X_3^* + s X_4^*:\ s, t \in \dsr}.
 \]
  On the
other hand, it is easy to verify that
\[
\frak g_\Gamma^* =\zspann\set{X_1^*, \ldots, X_3^*, \frac{1}{6}
X_4^*}.
\] Therefore
\[
\orbite_{\pi_l}^G\cap \frak g_\Gamma^* = \set{ X_1^*+  t X_2^* +
 \frac{t^2}{2}  X_3^* + \frac{s}{6} X_4^*:\ s\in \dsz, t \in
2\dsz}.
\]
 Let
\[
f_{t_0, s_0}= X_1^*+  t_0 X_2^* + \dfrac{t_0^2}{2} X_3^* +
\dfrac{s_0}{6} X_4^*\in \orbite_{\pi_l}^G\cap \frak g_\Gamma^*
\]
 and
\[
\gamma = \exp{rX_2}\exp{sX_3}\exp{6t X_4}\in \Gamma.
\]
We
calculate
\[
\Ad^*(\gamma) f_{t_0, s_0} = X_1^*+  (t_0-6t) X_2^* + \dfrac{(t_0-6t)^2}{2}
 X_3^*+
\left(\dfrac{s_0}{6}+st_0+r-6st\right) X_4^*.
\] Then (see \cite{Nielsen5})
\begin{gather*}
  \Ad^*(\Gamma) f_{t_0, s_0}  =  \set{ X_1^*+ (t_0+6t) X_2^* +
 \frac{(t_0+6t)^2}{2}  X_3^* + \left(\frac{s_0}{6}+ s\right)
X_4^*:\ s, t \in \dsz} \\
\phantom{\Ad^*(\Gamma) f_{t_0, s_0}}{} = \set{ f_{t_0+6t, s_0+6s}:\ s, t \in \dsz}.
\end{gather*} From this we deduce
that $\#[\orbite_{\pi_l}^G\cap \frak g_\Gamma^*/\Gamma] = 3\cdot 6 = 18$,
and hence
\begin{gather*}
\mult(\pi_l, G, \Gamma, 1)^2 \ne \#[\orbite_{\pi_l}^G\cap \frak
g_\Gamma^*/\Gamma].
\end{gather*}
Therefore, the group $G$ does not satisfy the Moore formula at~$\Gamma$.

  \subsection*{Acknowledgements}
   It is great pleasure  to thank the anonymous
  referees for their critical and valuable comments.

\pdfbookmark[1]{References}{ref}
\LastPageEnding

\end{document}